\documentclass{article}
\usepackage{amsmath}
\usepackage{amssymb}
\newtheorem{thm}{THEOREM}[section]
\newtheorem{lem}[thm]{LEMMA}
\newtheorem{prop}[thm]{PROPOSITION}
\newtheorem{cor}[thm]{COROLLARY}

\begin{document}
\begin{flushleft}
\huge{\bf ON PLANE ALGEBROID CURVES}
\end{flushleft}

\bigskip\noindent
V. BARUCCI, Universit\`a di Roma 1,
{\em email}: barucci@mat.uniroma1.it

\medskip\noindent
M. D'ANNA, Universit\`a di Catania,
{\em email}: mdanna@dmi.unict.it

\medskip\noindent
R. FR\"OBERG,
Stockholms Universitet,
{\em email}: ralff@matematik.su.se

\begin{abstract}
Two plane analytic branches are topologically equivalent if and only
if they
have the same multiplicity sequence. We show that having same semigroup
is equivalent to having same multiplicity sequence, we calculate the
semigroup from a parametrization, and we characterize semigroups for plane
branches. These results are known, but the proofs are new. 
Furthermore we characterize multiplicity sequences of plane branches,
and we prove that the associated graded ring, with respect to the values,
of a plane branch is a complete intersection.
\end{abstract}

\section{INTRODUCTION}

Let $C$ and $C'$ be two analytic plane irreducible curves (branches)
defined in a neighbourhood of the origin and having singularities
there. The branches are said to be topologically equivalent
if there are neighbourhoods $U$ and $U'$ of the origin
such that $C$ is defined in $U$, $C'$ in $U'$, and there is a homeomorphism
$T\colon U\rightarrow U'$ such that $T(C\cap U)=C'\cap U'$.

If $F(X,Y)\in{\mathbb C}[[X,Y]]$ is an irreducible formal power series,
the local ring ${\mathcal O}={\mathbb C}[[X,Y]]/(F)$ is called a (plane)
algebroid branch. Two algebroid branches are formally equivalent if
they have the same multiplicity sequence (see below for the definition
of multiplicity sequence). Every algebroid (analytic resp.) branch is
 formally (topologically, resp.) equivalent to an algebraic branch,
i.e. a branch defined by a polynomial \cite{Sa}, and if two analytic
branches are formally equivalent, they are topologically equivalent. We will
in the sequel consider algebroid branches.

Zariski has shown (\cite{Za}) that two branches are formally equivalent
if and only if they have the same semigroup of values (see below for
the definition of the value semigroup of a branch).

The crucial result of Section 2 is Proposition \ref{ap}, which gives the
relation between the value semigroups of an algebroid plane branch
${\mathcal O}$ and its blowup ${\mathcal O'}$. It is a result contained in
\cite{Ap}. Ap\'ery proved that, in order to show that
the value semigroup $v({\mathcal O})$ of an algebroid plane branch
${\mathcal O}$ is  symmetric. Subsequently Kunz proved that, for any
analytically irreducible ring
${\mathcal O}$, ${\mathcal O}$ is Gorenstein if and only if $v({\mathcal
O})$ is symmetric. So now it is more common to say that the value
semigroup of an algebroid plane branch is symmetric because the ring is
Gorenstein (it is in fact a complete intersection). At any rate we are
interested in Ap\'ery's result for different reasons. By its use we give an
easy proof of the fact that two plane algebroid
branches are formally equivalent if and only if they have the same
semigroup of values. We get also a well known formula of Hironaka and
apply it again in Sections 3 and 4.  The material in Section 3 is
classical too and essentially contained in Enriqes-Chisini's work, but what
is new, is the use of Ap\'ery's Lemma in this context. After characterizing
all possible multiplicity sequences for plane branches, we give a
criterion to check if a semigroup is the value semigroup of a plane branch.
In Section 4, we determine the semigroup of a plane branch from its
parametrization,  here  also using results from
\cite{Ap}. This result is well known, but the proof is new as far as we
know. Finally in Section 5 we show that the semigroup ring of the semigroup
of a plane curve is a complete intersection.

\section{PLANE BRANCHES}

Starting from Ap\'ery's article \cite{Ap}, we will proceed to explicate
and expand various elements that are presented in the original arguments
in a summary or not totally developed manner.

\medskip
Let ${\mathcal O}=
{\mathbb C}[[X,Y]]/(F)={\mathbb C}[[x,y]]$, where $F$ is irreducible in
${\mathbb C}[[X,Y]]$ be an algebroid plane branch.
Since $F(X,Y)$ is irreducible, then
$F(X,Y)$ must contain some term $X^i$ and some term $Y^j$
(otherwise $F$ is not irreducible since we could factor out $X$
or $Y$). Denote the minimal such powers by $n$ and $m$ respectively.
Then, by the
Weierstrass Preparation Theorem,  the same ideal $(F)$ can be generated by
an element $X^n+\phi(X,Y)$, where $\phi(X,Y)$ is a polynomial
of degree $n-1$ in $X$ with coefficients which are power series in
$Y$ (or vice versa by an element $Y^m+\psi(X,Y)$, where $\psi(X,Y)$
is a polynomial of degree $m-1$ in $Y$ with coefficients which are
power series in $X$). This gives  that ${\mathcal O}$ is generated by
$1,x,...,x^{n-1}$ as ${\mathbb C}[[y]]$-module (or generated by
$1,y,...,y^{m-1}$ as ${\mathbb C}[[x]]$-module).

The Puiseux Theorem
gives that the branch has a parametric representation $x=t^m,y=\sum a_it^i$
(or $x=\sum_{i \geq m} b_it_1^i,y=t_1^n$,
where ${\mathbb C}[[t]]={\mathbb C}[[t_1]]$).
Thus ${\mathcal O}={\mathbb C}[[x,y]] \subseteq
{\mathbb C}[[t]] = \bar{\mathcal O}$,
which is a discrete valuation
ring. Denote by $v$ the valuation of such ring that consists in
associating to any formal
power series in ${\mathbb C}[[t]]$ its order.
In particular $v(x)=m$ and $v(y)=n$.
Since the fraction field of
${\mathcal O}$ equals the fraction field of
$\bar{\mathcal O}$, there exist $f_1(t),f_2(t)\in{\mathcal O}$, such
that $f_1(t)/f_2(t)=t$, so $f_1(t)=tf_2(t)$ and $v(f_1)=v(f_2)+1$.
Since $\gcd(v(f_1),v(f_2))=1$, all sufficiently large integers belong
to $v({\mathcal O})=\{ v(z);z\in{\mathcal O}\setminus\{ 0\}\}$. Thus
$v({\mathcal O})$ is a numerical semigroup, i.e., a subsemigroup of
${\mathbb N}$ with finite complement to ${\mathbb N}$.

In the sequel we use the following terminology.  If $S$ is a subsemigroup
of $\mathbb N$ and $T$ is a subset of
$\mathbb Z$, we call $T$ an $S$-module if $s\in S,t\in T$ implies
$s+t\in T$. We call $T$ a free $S$-module if $T=\cup_{i=1}^kT_i$
with $T_i\cap T_j=\emptyset$ if $i\ne j$
and $T_i=n_i+S$ for some $n_i\in{\mathbb Z}$. We call $n_1,\ldots,n_k$
a basis of $T$.

\medskip
With the hypotheses and notation above, we will construct a new basis
$y_0,\ldots,y_{m-1}$ for
${\mathcal O}$ as a
${\mathbb C}[[x]]$-module, such that, for each $i$, $y_0,\ldots,y_i$
is a basis for ${\mathcal O}_i={\mathbb C}[[x]]+y{\mathbb C}[[x]]+
\cdots+y^i{\mathbb C}[[x]]$, and furthermore such that $v({\mathcal O}_i)=
\{ v(z);z\in{\mathcal O}_i\setminus\{ 0\}\}$ is a free module over
$v({\mathbb C}[[x]])=m{\mathbb N}$ with basis $\omega_0,\ldots,\omega_i$,
where each $\omega_j=v(y_j),\, j=0,\ldots,i$ is the smallest value
in $v({\mathcal O})$ in its congruence class $\pmod{m}$.
Let $y_0=1$, thus $\omega_0=v(y_0)=0$ and $v({\mathcal O}_0)=
v({\mathbb C}[[x]])=m{\mathbb N}$. Suppose that
$y_0,\ldots,y_{k-1},k<m$ have been defined such that $v({\mathcal O}_{k-1})$
is
a free $m{\mathbb N}$-module with basis $\omega_0,\ldots,\omega_{k-1}$.
 We claim that there exists a
$\phi(x,y)\in{\mathcal O}_{k-1}$ such that $y_k=y^k+\phi(x,y)$
has a value which does not belong to
$v({\mathcal O}_{k-1})$. If $v(y^k)\notin v({\mathcal O}_{k-1})$,
we are ready. Otherwise
$v(y^k)=v(z_1)$ for some $z_1\in{\mathcal O}_{k-1}$.
Then $v(y^k-c_1z_1)>v(y^k)$ for some
$c_1\in{\mathbb C}$. If $v(y^k-c_1z_1)\notin v({\mathcal O}_{k-1})$,
we are ready. Otherwise
take $z_2\in{\mathcal O}_{k-1}$ with $v(z_2)=v(y^{k}-c_1z_1)$. Then
$v(y^k-c_1z_1-c_2z_2)>v(y^k-c_1z_1)$ for some $c_2\in{\mathbb C}$ a.s.o.
Thus we see that the expansion of $y^k$ as a power series in $t$ must
contain a term $a_it^i$ with $a_i\ne0$ and $i\notin v({\mathcal O}_{k-1})$,
since otherwise $y^k\in{\mathcal O}_{k-1}$.

Notice that  $y_1y_{k-1}=(y+\phi_1(x))(y^{k-1}+\phi_{k-1}(x,y))=
y^k+\psi(x,y),\psi(x,y)\in{\mathcal O}_{k-1}$, so $y_k =
y_1y_{k-1}+\phi(x,y) - \psi(x,y)$
and we could equally well
have defined $y_k$ as an element of the form $y_1y_{k-1}+\phi(x,y)$
(where $\phi(x,y)
\in{\mathcal O}_{k-1}$) with a value which does not belong to
$v({\mathcal O}_{k-1})$. In such expression of $y_k$,  $v(\phi(x,y)) \geq
v(y_1y_{k-1})$
 since otherwise $v(y_k)=v(\phi(x,y))
\in v({\mathcal O}_{k-1})$. Thus
$\omega_k=v(y_k)\ge v(y_1y_{k-1})=v(y_1)+v(y_{k-1})=
\omega_1+\omega_{k-1}$. In particular the sequence
$\omega_0,\omega_1,\ldots$
is strictly increasing. Since $v({\mathcal O}_{k-1})$ is free over
$m{\mathbb Z}$, this shows that $\omega_k\not\equiv\omega_j$ if
$j<k$.
 Any element $z\in{\mathcal O}_k$ can be written
$z=a_0(x)y_0+\cdots+a_k(x)y_k$. All terms in this sum have values in
different congruence classes $\pmod{m}$. Thus $v(z)=\min v(a_i(x)y_i)$.
This shows that $v({\mathcal O}_k)$ is free with basis $\omega_0,\ldots,
\omega_k$. After $m$ steps, we get that ${\mathcal O}_{m-1}$=${\mathcal O}$
is a
${\mathbb C}[[x]]$-module generated by $y_0,\ldots,y_{m-1}$ with the
requested properties.

\medskip
If $S$ is a numerical semigroup and
$a\in S\setminus\{ 0\}$, then the
elements $n_0,n_1,\ldots,\\ n_{a-1}$, where $n_i$ is the smallest element
in $S$ congruent to $i\pmod a$, is called the {\it Apery set} of $S$ with
respect to $a$.
If we order the elements in the Apery set, and then denote them $\omega_0,
\ldots,\omega_{a-1}$, we have the {\em ordered Apery set}.
We call the elements $y_0,\ldots,y_{m-1}\in{\mathcal
O}$ constructed as above
an {\it Apery basis} of ${\mathcal O}$ with respect to $x$ .
By the construction,
$\omega_0=v(y_0),
\ldots,\omega_{m-1}=v(y_{m-1})$ is the ordered Apery set of $v({\mathcal
O})$.

In a similar way an Apery basis of ${\mathcal O}$ with respect to $y$ is
defined.

\medskip\noindent{\bf Example} If in ${\mathcal O}={\mathbb C}[[x,
y]]$ we have $\gcd(m,n)=1$, where $v(x)=m$ and
$v(y)=n$,
then $y_k=y^k,\, k=0,\ldots,m-1$ is an Apery basis of ${\mathcal O}$, and
thus $\omega_k=kn,\, k=0,\ldots,m-1$ is the ordered Apery set of
$v({\mathcal O})$ with respect to $m$.

\medskip\noindent{\bf Example} If in ${\mathcal O}={\mathbb C}[[x,
y]]$ we have $x=t^8,\, y=t^{12}+t^{14}+t^{15}$,
then
$y_0=1,\, y_1=y,\, y_2=y^2-x^3=2t^{26}+\cdots,\, y_3=y^3-x^3y=2t^{38},\,
y_4=y^4-2x^3y^2-4x^5y+3x^6=8t^{53}+\cdots,\, y_5=y^5-2x^3y^3+x^6y-4x^8=
8t^{65}+\cdots,\,
y_6=y^6-3x^3y^4-4x^5y^3+3x^6y^2+4x^8y-x^9=16t^{79}+\cdots$,
and $y_7=y^7-3x^3y^5+3x^6y^3-4x^8y^2-x^9y+4x^{11}=16t^{91}
+\cdots$ is an Apery basis for ${\mathcal O}$, so the ordered
Apery set of $v({\mathcal O})$ with respect to $8$
is $\{ 0,12,26,38,53,65,79,91\}$. Thus $v({\mathcal
O})$ is minimally generated by $8,12,26,53$.

\medskip
If $S$ is a numerical semigroup, we denote the {\it Frobenius number} of $S$,
i.e. $\max\{ x\in{\mathbb Z};\, x\notin S\}$, by $\gamma(S)$.
The {\it conductor}
of $S$ is $c(S)=\gamma(S)+1=\min\{ x;\, [x,\infty)\subseteq S\}$.

The following lemma is well known, and its easy proof is left to the reader.

\begin{lem}\label{Frobenius}
Let $S$ be a numerical semigroup with Frobenius number $\gamma$ and $a\in
S$.
If $\omega_0,\ldots,\omega_{a-1}$ is the ordered Apery set of $S$ with
respect to $a$, then $\gamma=\omega_{a-1}-a$.
\end{lem}

Now we are ready for the crucial lemma from \cite{Ap}.
If ${\mathcal O}={\mathbb C}[[x,
y]]$ with $v(x)<v(y)$, we denote the quadratic transform (or blowup)
${\mathbb C}[[x,y/x]]$ by ${\mathcal O}'$.

\begin{lem}
If an Apery basis of ${\mathcal O}'$ with respect to $x$ is
$y_0',\ldots,y_{m-1}'$, then $y_i=y_i'x^i$, for $i=0,\ldots,m-1$ is
an Apery basis of ${\mathcal O}$ with respect to $x$.
\end{lem}

\noindent{\bf Proof.} Let $F_i(x,y/x)$ be the polynomial of degree $i$ in
$y/x$ which defines $y'_i$, i.e. let $F_i(x,y/x) = (y/x)^i+ \phi'_i(x,y/x)$,
where $\deg(\phi'_i)<i$ in $y/x$.
Then $y_i=x^iF(x,y/x)=y^i+\phi_i(x,y),
\phi_i(x,y)\in{\mathcal O}_{i-1}$ is of the requested form and,
if $v(y_i')=\omega_i'$, then $\omega_i=v(y_i)=\omega_i'+im$, thus
$\omega_i\equiv\omega_i'\pmod m$. We have to show that $\omega_i\notin
v({\mathcal O}_{i-1})$. This is because $\omega_i'$ is not congruent
to any $\omega_j'$, if $j<i$, and so also $\omega_i$ is not congruent
to any $\omega_j$, if $j<i$.

\medskip
As a consequence we get

\begin{prop}\label{ap}\cite[Lemme 2]{Ap}
If the ordered Apery set of $v({\mathcal O}')$
with respect to $m=v(x)$
is $0=\omega'_0<\omega'_1<\cdots<\omega'_{m-1}$, then the ordered
Apery set of
$v({\mathcal O})$ with respect to $m$ is $\omega_0=\omega'_0<\omega_1=
\omega'_1+m<\omega_2=\omega'_2+2m<\ldots<\omega_{m-1}=\omega'_{m-1}+(m-1)m$.
\end{prop}

Recall that the multiplicity of the ring ${\mathcal O} =
 {\mathbb C}[[x,y]]$, where  $x=a_mt^m+a_{m+1}t^{m+1}+\cdots,\, a_m\ne0$
and $y=b_nt^n+b_{n+1}t^{n+1}
+\cdots,\, b_n\ne0$, is given by $\min(m,n)$ i.e. the multiplicity of
${\mathcal O}$
 is the smallest positive value in $v({\mathcal O})$.

Set ${\mathcal O}={\mathcal O}^{(0)}$, denote by ${\mathcal O}^{(i+1)}$ the
blowup of
 ${\mathcal O}^{(i)}$ and by $e_i$ the multiplicity of ${\mathcal O}^{(i)}$.
The {\it multiplicity
sequence} of ${\mathcal O}$ is by definition the sequence of natural numbers
$e_0, e_1, e_2, \cdots$.
Let $k$ be the minimal index such that $e_k = 1$, i.e. such that
$v({\mathcal O}^{(k)})= {\mathbb N}$. Two algebroid branches are {\it
formally equivalent} if they have the same multiplicity sequence.

As a consequence of Proposition \ref{ap}, we get easily a well known
formula:

\begin{cor}\label{hironaka}\cite[Theorem 1]{Hi}
We have $l_{\mathcal O}(\bar{\mathcal O}/({\mathcal O}:\bar{\mathcal O}))
=\sum_{i=0}^{k} e_i(e_i-1)$ and
$l_{\mathcal O}({\mathcal O}/({\mathcal O}:\bar{\mathcal O}))=
l_{\mathcal O}(\bar{\mathcal O}/{\mathcal O})=\frac{1}{2}\sum_{i=0}^{k}
e_i(e_i-1)$.
\end{cor}

\noindent{\bf Proof.} Let $\omega_i^{(j)}$ ($\omega_i^{(j+1)}$, resp.)
be the $i$'th element in the ordered Apery set of
$v({\mathcal O}^{(j)})$, ($v({\mathcal O}^{(j+1)})$, resp.),  with respect
to $e_j$  and let
${\mathcal O}^{(j)}:\bar{\mathcal O}=t^{c_j}{\mathbb C}[[t]]$
(${\mathcal O}^{(j+1)}:\bar{\mathcal O}=t^{c_{j+1}}{\mathbb C}[[t]]$,
resp.). By Lemma \ref{Frobenius}  $c_j=\omega_{e_j-1}^{(j)}-e_j+1$ and
$c_{j+1}=\omega_{e_j-1}^{(j+1)}-e_j+1$. Proposition \ref{ap} gives
$\omega_{e_j-1}^{(j)} = \omega_{e_j-1}^{(j+1)}+e_j(e_j-1)$ and so
$c_j=c_{j+1}+e_j(e_j-1)$. It follows that $c_0=l_{\mathcal O}(\bar{\mathcal
O}/({\mathcal O}:\bar{\mathcal O}))
=c_1+e_0(e_0-1)= \cdots=c_k+e_{k-1}(e_{k-1}-1)+
\cdots+e_0(e_0-1)=\sum_{i=0}^{k} e_i(e_i-1)$. Since the ring ${\mathcal O}$
is
Gorenstein, we get $l_{\mathcal O}({\mathcal O}/({\mathcal O}:\bar{\mathcal
O}))=
l_{\mathcal O}(\bar{\mathcal O}/{\mathcal O})=\frac{1}{2}l_{\mathcal
O}(\bar{\mathcal O}/({\mathcal O}:\bar{\mathcal O}))$.

\medskip
\noindent{\bf Example} Not every symmetric semigroup is the value
semigroup of an algebroid plane branch. The semigroup generated by 4,5,6 is
symmetric and has Apery set 0,5,6,11 with respect to 4. If this were the
value semigroup of a plane branch, then the Apery set of its blowup would
be
$0,\, 1=5-4,\, -2=6-8,\, -1=11-12$ which obviously is impossible.

\begin{thm}\label{a-equiv} \cite{Za} Two algebroid plane branches are
formally equivalent if and only if they have the same semigroup.
\end{thm}

\noindent{\bf Proof.} Let ${\mathcal O}={\mathcal O}^{(0)},
{\mathcal O}^{(1)},\ldots$ be the sequence of blowups of $\mathcal O$,
and let $e_0,\ldots,e_k=1$ be the corresponding multiplicity sequence.
Then $v({\mathcal O}^{(k)})=
{\mathbb N}$ has ordered Apery set $\{0,1,\ldots,e_{k-1}-1\}$
with respect to $e_{k-1}$. Proposition \ref{ap} gives the ordered Apery
set, hence the semigroup, of ${\mathcal O}^{(k-1)}$ with respect to
$e_{k-1}$ a.s.o. Thus the multiplicity sequence determines the semigroup of
${\mathcal O}$. On the other hand, the semigroup of ${\mathcal O}$
gives the multiplicity $e_0$ of ${\mathcal O}$. Proposition \ref{ap} gives the
Apery set of $v({\mathcal O}^{(1)})$, hence $v({\mathcal O}^{(1)})$
and so on. Thus the semigroup $v({\mathcal O})$ gives the
multiplicity sequence.

\medskip
Let $c_i$ denote the {\it conductor
degree} of ${\mathcal O}^{(i)}$, i.e. ${\mathcal O}^{(i)}:\bar{\mathcal O}=
t^{c_i}{\mathbb C}[[t]]$, and call $(c_0,c_1,\ldots)$ the {\it conductor
degree sequence} of $\mathcal O$. Let $f_i=l_{{\mathcal O}^{(i)}}(
\bar{\mathcal O}/{\mathcal O}^{(i)})$, and call $(f_0,f_1,\ldots)$
the {\it sequence of singularity degrees} of $\mathcal O$.

\begin{cor} Two algebroid plane branches
are
formally equivalent if and only if they have the same conductor degree
sequence, and if and only if they have the same sequence of singularity
degrees.
\end{cor}

\noindent{\bf Proof.} If $\omega_0<\omega_1<\cdots<\omega_{e_i-1}$ is the
Apery set of $v({\mathcal O}^{(i)})$ with respect to $e_i$, then by
Lemma \ref{Frobenius} $c_i=
\omega_{e_i-1}-e_i+1$. Thus the multiplicity sequence of $\mathcal O$
determines, and is determined by, the conductor degree sequence. Since each
ring ${\mathcal O}^{(i)}$ is Gorenstein,
$f_i=c_i/2$ and the same is true for the sequence of singularity degrees.

\medskip\noindent{\bf Example} The conductor degree of $\mathcal O$ does
not suffice to give formal equivalence. The branches ${\mathbb C}[[t^4,t^5]]$
and ${\mathbb C}[[t^3,t^7]]$ both have conductor $t^{12}{\mathbb C}[[t]]$,
but they are not formally equivalent.

\section{THE MULTIPLICITY SEQUENCE FOR A\\ PLANE BRANCH}

A sequence of numbers $e_0\ge e_1\ge e_2\ge\cdots$ is a
multiplicity sequence of a (not necessarily plane) branch if and only if
$0,e_0,e_0+e_1,e_0+e_1+e_2,\ldots$ constitute a semigroup \cite{dV}. We will
now determine which multiplicity sequences occur for plane branches.
We will also use this result together with Proposition \ref{ap} and
Theorem \ref{a-equiv}
to get an algorithm to determine if a symmetric semigroup is
the semigroup of a plane branch.

Let ${\mathcal O}=
{\mathbb C}[[t^{\delta_0}+\sum_{i\ge N}a_it^i,
\sum_{i>\delta_0}b_it^i]]$ be a branch.
Let, for $i\ge1$,
$\delta_i=\min\{ j;\, b_j\ne0,\gcd(\delta_0,\ldots,\delta_{i-1},j)<
\gcd(\delta_0,\ldots,\delta_{i-1})\}$.
Let $d_0=\delta_0$ and $\gcd(\delta_0,\ldots,\delta_i)=d_i$ for $i\ge1$. Set
 also $k=\min
\{ i;\, d_i=1\}$. (There exists such a $k$ since the integral
closure of ${\mathcal O}$ is ${\mathbb C}[[t]]$.)
We call the parametrization {\em standard} if $N>\delta_k$.
The numbers $\delta_0,\delta_1,\ldots$ are called the {\it characteristic
exponents} of $\mathcal O$. It follows from the proof of Lemma \ref{newpar}
below, that we always can get a standard parametrization from a given one.

\begin{lem}\label{newpar}
Let ${\mathcal O}=
{\mathbb C}[[t^{\delta_0}+\sum_{i\ge N}a_it^i,
\sum_{i>\delta_0}b_it^i]]$ be a branch with standard parametrization and
with characteristic exponents $(\delta_0,\ldots,\delta_k)$.
Then the characteristic exponents of $\mathcal O'$ are:\\
a) $(\delta_0,\delta_1-\delta_0,
\ldots,\delta_k-\delta_0)$, if $\delta_0<\delta_1-\delta_0$.\\
b) $(\delta_1-\delta_0,\delta_0,\delta_0+\delta_2-\delta_1,
\ldots,\delta_0+\delta_k-\delta_1)$, if $\delta_0>\delta_1-\delta_0$
and $\delta_0$ is not a multiple of
$\delta_1-\delta_0$\\
c) $(\delta_1-\delta_0,\delta_0+\delta_2-\delta_1,\ldots,\delta_0+\delta_k-
\delta_1)$, if $\delta_0$ is a multiple of $\delta_1-\delta_0$.\\
\end{lem}

\noindent{\bf Proof}
We can suppose that
$v(\sum_{i>\delta_0}b_it^i)=\delta_1$. Then the blowup
${\mathcal O}'$ of ${\mathcal O}$ is
$(t^{\delta_0}+\cdots,t^{\delta_1-\delta_0}+\cdots)$.
One of the following three cases will occur:\\
a) $\delta_0<\delta_1-\delta_0$\\
b) $\delta_0>\delta_1-\delta_0$ and $\delta_0$ is not a multiple of
$\delta_1-\delta_0$\\
c) $\delta_0$ is a multiple of $\delta_1-\delta_0$.\\
We will in each case write ${\mathcal O}'$ in standard form
and derive its characteristic exponents.
In case a) ${\mathcal O}'$ is of standard form. We keep the meaning of
$\delta_i$ and $d_i$ from above and denote the
corresponding entities for ${\mathcal O}'$ with $\delta_i'$ and $d_i'$.
It follows that $d_i'=d_i$ for all $i$ and that ${\mathcal O}'$ has
characteristic exponents
$(\delta_0',\delta_1',\ldots,\delta_k')=(\delta_0,\delta_1-\delta_0,
\ldots,\delta_k-\delta_0)$. In case b) we first make the coordinate change
$X=y,Y=x$ to get $(t^{\delta_1-\delta_0}(1+
\sum_{i\ge1}c_it^i),t^{\delta_0}+\cdots)$. Let $i_0=\min\{ i;\, c_i\ne0\}$.
Then we choose a new parameter $t_1$, by $t=t_1(1-\frac{c_{i_0}}
{\delta_1-\delta_0}t_1^{c_{i_0}})$ to get the parametrization
$(t_1^{\delta_1-\delta_0}(1+\sum_{i\ge1}c_i't_1^i),t_1^{\delta_0}+\cdots)$.
Now $v(\sum_{i\ge1}c_i't_1^i)>v(\sum_{i\ge1}c_it^i)$.
We continue to change parameter in this way. After a finite number of
steps we get a parametrization of the branch of the type
$(t^{\delta_1-\delta_0}+
\sum_{i>\delta_k}k_it^i,t^{\delta_0}+\cdots)$ with $d_i'=d_i$ for
all $i$, and with characteristic exponents
$(\delta_1-\delta_0,\delta_0,\delta_0+\delta_2-\delta_1,
\ldots,\delta_0+\delta_k-\delta_1)$. In case c) finally, we use a similar
reparametrization and get $d_i'=d_{i+1}$, and a branch with characteristic
exponents
$(\delta_1-\delta_0,\delta_0+\delta_2-\delta_1,\ldots,\delta_0+\delta_k-
\delta_1)$.

\medskip     

If $m_0, m_1, \cdots$ and $h_0, h_1, \cdots$ are natural numbers, denote by
$m_0^{(h_0)}, m_1^{(h_1)}, \cdots$ the sequence of
natural numbers given by $m_0$ repeated $h_0$ times, $m_1$ repeated
$h_1$ times and so on.
Suppose that for a couple $m,n$ of natural numbers,
the Euclidean algorithm gives

$$ m=nq_1+r_1$$
$$ n=r_1q_2+r_2$$
$$\cdots$$
$$ r_{i-1}=r_iq_{i+1}+r_{i+1}$$
$$ r_i=r_{i+1}q_{i+2}+0$$

Denote by $M(m,n)$ the sequence of natural numbers $n^{(q_1)}, r_1^{(q_2)},
\cdots, r_{i+1}^{(q_{i+2})}$. Of course such a sequence
ends with $r_{i+1}= \gcd(m,n)$ (if $m<n$, and so $q_1=0$, $n$ appears $0$
times, i.e. it does not appear, hence
$M(m,n)=M(n,m)$). With this notation:

\begin{thm}\label{mseq} A sequence of natural numbers is the
multiplicity sequence  of
an algebroid plane branch if and only if it is of
the following form:
$$M(m_0,m_1), M(m_2,m_3), \ldots,M(m_{2k},m_{2k+1}),1,1,\ldots $$
where, for $i \geq 0$,  $\gcd (m_{2i},m_{2i+1})=m_{2i+2}$ and $m_{2i+3}$ is
such that $m_{2i+4} < m_{2i+2}$, and finally $\gcd(m_{2k},m_{2k+1})=1$.
\end{thm}

\noindent{\bf Proof.} Let ${\mathcal O}$ be an algebroid plane branch
with standard parametrization. Then, by Lemma \ref{newpar}, its
multiplicity sequence is
$$M(\delta_0, \delta_1), M(d_1, \delta_2 - \delta_1), M(d_2, \delta_3 -
\delta_2),\ldots,M(d_{k-1},\delta_k-\delta_{k-1}),1,1,\ldots$$
and is a sequence of the requested form. Conversely, given a sequence of
natural numbers as in the statement, we can get
characteristic exponents
$(\delta_0, \delta_1, \ldots \delta_k)$ and so an ${\mathcal O}$.

\medskip
We give two concrete examples.

\medskip
\noindent{\bf Example} $6,4,2,2,1,1,\ldots=M(10,6),1,1,\ldots$ is an
admissible multiplicity sequence (i.e. the multiplicity sequence of an
algebroid plane branch), but
$6,4,2,1,1,\ldots$ is not.

\medskip
\noindent{\bf Example} Let
$\mathcal O={\mathbb C}[[x,y]]$ with
$$x=t^{2\cdot 2^n},\, y=t^{3\cdot 2^n}+t^{3\cdot 2^n+2^{n-1}}+\cdots+
t^{3\cdot 2^n+2^{n-1}+\cdots+2+1}.$$
The multiplicity sequence is
$$2^{n+1},2^n,2^n,2^{n-1},2^{n-1},\ldots,4,4,2,2,1,\ldots=$$
$$M(3\cdot2^n,2^{n+1}),M(2^n,2^{n-1}),M(2^{n-1},2^{n-2}),\ldots,M(2,1),\ldots$$

\medskip

Now we are ready to give an algorithm to determine if
a symmetric semigroup is the semigroup of values of a plane curve.

\begin{lem}\label{step}
Let $S$ be a symmetric semigroup, $m=\min(S\setminus \{0\})$ and
let
$$0=\omega_0<\omega_1<\cdots<\omega_{m-1}$$
be its ordered Apery set with respect to $m$.
Suppose that $\omega_0<\omega_1-m<\cdots<\omega_{m-1}-(m-1)m$
is the ordered Apery set of a semigroup $S'$.
Then $S'$ is symmetric.
\end{lem}

\noindent{\bf Proof.}
This follows from \cite{Ap}. 

\medskip
Given a symmetric semigroup $S$ satisfying the hypotheses of
Lemma \ref{step}, one could repeat the process for the ordered
Apery set of $S'$ with respect to its minimal non zero element, and so on
until we find either a semigroup which does not satify these hypotheses
or we find $\mathbb N$.
But even if, after a finite number of steps, we get $\mathbb N$,
it is not true that $S$ is a value semigroup of a plane branch,
as the following example shows.

\medskip
\noindent{\bf Example}
Let $S=\langle6,10,29\rangle$; its ordered Apery set with respect to
$6$ is $\{0,10,20,\\ 29,39,49\}$. The set  obtained applying
Lemma \ref{step} is $\{0,4=10-6,8=20-12,11=29-18,15=39-24,19=49-30\}$,
hence it is the ordered Apery
set of $S'$ with respect to $6$.
Hence
$S'=\{0,4,6,8,10,11,12,14,\rightarrow \dots\}=\langle4,6,11\rangle$.
The ordered Apery set of $S'$ with respect to $4$ is
$0,6,11,17$. Hence we get the new set $\{0,2=6-4,3=11-8,5=17-12\}$
which is still ordered
and determines the semigroup $S''=\langle2,3\rangle$. Its ordered
Apery set with respect to $2$ is $0,3$. Thus we get the set
$\{0,1=3-2\}$, which is the ordered Apery set with respect to
$2$ of $\mathbb N$.

On the other hand the semigroup $S$ is not the value
semigroup of a plane branch $\cal O$
since the multiplicity sequence of $\cal O$ should be $6,4,2,1,1,\dots$
which is not admissible, since the subsequence $6,4$
can be obtained only by $M(10,6)$ but $M(10,6)= 6,4,2,2$.

\bigskip
Let $S$ be the value semigroup of a plane branch $\mathcal O$. By Proposition
\ref{ap} we get that $S'$ (defined as in Lemma \ref{step}) is again a
symmetric semigroup
and $S'=v({\mathcal O}')$. Repeating the process, if $S^{(0)}=S$
and $S^{(j+1)}=(S^{(j)})'$, and denoting by $m_j$ the minimal non zero element
of $S^{(j)}$ and by $\omega_0^{(j)},\omega_1^{(j)},
\ldots,\omega_{m_j-1}^{(j)}$ its
ordered Apery set with respect to $m_j$, we get that
$\omega_0^{(j)},\omega_1^{(j)}-m,\ldots,\omega_{m_j-1}^{(j)}-(m_j-1)m_j$ is
the ordered Apery set of a symmetric semigroup
$S^{(j+1)}$, and $S^{(j+1)}=
v({\mathcal O}^{(j+1)})$. Since there exists an $n\ge1$ such that
${\mathcal O}^{(n)}={\mathbb C}[[t]]$, then $S^{(n)}={\mathbb N}$. Moreover
the sequence $m_0,\ldots,m_{n-1},1,\ldots$ is the multiplicity sequence
of ${\mathcal O}$, hence is an admissible multiplicity sequence.

Conversely if $S=S_0$ is a symmetric semigroup, let $S^{(j)},m_j,
\omega_i^{(j)}$ be defined as above.
If the sets $0=\omega_0^{(j)},
\omega_1^{(j)}-m_j, \dots, \omega_{m_j-1}^{(j)}-(m_j-1)m_j$
are ordered Apery sets for every $j=0,\dots n-1$ and the
sequence $m_0, m_1, \dots, m_{n-1}, 1,1,\dots$ is an admissible
multiplicity sequence, then $S$ is the value semigroup of a plane branch.
In fact,
since the sequence
$m_0,m_1, \dots, m_{n-1},1,1,\ldots$ is an admissible multiplicity
sequence, then there exists a plane branch $\mathcal O$
having this sequence as multiplicity sequence.
Now, by
Theorem \ref{a-equiv}, the multiplicity sequence determines
the value semigroups $v(\mathcal O^{(k)})$, $k=0\dots,n-1$,
and these semigroups, by Proposition \ref{ap} and Lemma \ref{step}
have the same ordered  Apery sets of the semigroups $S^{(k)}$;
hence they are the same semigroups.

This discussion gives a criterion to check if $S$ is the value semigroup
of a plane branch, since we can apply repeatedly the process described in
Lemma \ref{step} until we find either a semigroup which does not
satisfy the hypotheses in Lemma \ref{step}
or we find ${\mathbb N}$. If the last case occurs, then it is enough to
check if the sequence $m_0,\ldots,m_{n-1},1,1,\ldots$ is admissible.

\medskip
The condition that at each step
the sequence $0=\omega_0^{(j)},
\omega_1^{(j)}-m_j, \dots, \omega_{m_j-1}^{(j)}-(m_j-1)m_j$
is an ordered Apery set (and not only an Apery set) is necessary
as the following example shows.

\medskip
\noindent{\bf Example} Let $S=\{0,4,8,9,10,12,13,14,16,\rightarrow\dots\}$
be the semigroup with
ordered Apery set $\{0,9,10,19\}$ with respect to $4$.
The sequence $0,5=9-4,2=10-8,7=19-12$ is not increasing.
If we consider the semigroup $S'$ with ordered Apery set
$\{0,2,5,7\}$ with respect to $4$
it is the symmetric semigroup $\{0,2,4,\rightarrow \dots\}$
and then in two more steps we get $\mathbb N$.

Notice that the sequence $m_0, m_1, \dots$
is in this case $4,2,2,1,1,\dots$; it is admissible as multiplicity sequence
since it is
$M(6,4),M(2,1),1,1,\ldots$.
However, applying Theorem \ref{a-equiv}, we get the semigroup
$\{0,4,6,8,10,12,13,14,\rightarrow\dots\}$
with ordered Apery set $\{0,6,13,19\}$ and
applying Theorem \ref{mseq} we get the parametrization
$\mathcal O=\mathbb C[[t^4, t^6+t^7]]$.

\section{THE SEMIGROUP OF VALUES FOR A\\ PLANE BRANCH}

The following theorem is proved in different ways in e.g. \cite{Za},
\cite{Ca}, \cite{Ab},
\cite{Be-Ca}, \cite{Me}.

\begin{thm}\label{many} Let ${\mathcal O}=
{\mathbb C}[[t^{\delta_0}+\sum_{i\ge N}a_it^i,
\sum_{i>\delta_0}b_it^i]]$ be a branch with standard parametrization.
Denote the minimal generators of $v({\mathcal O})$ by
$\bar\delta_0<\cdots<\bar\delta_s$. Then $s=k$,
$\bar\delta_0=\delta_0,
\bar\delta_1=\delta_1$ and $\bar\delta_i=
\bar\delta_{i-1}\frac{d_{i-2}}{d_{i-1}}
+\delta_i-\delta_{i-1}$ if $i=2,\ldots,k$.
\end{thm}

We will divide the proof into several steps. From now on we will, for
a plane branch with characteristic exponents $(\delta_0,\delta_1,\ldots)$,
let $\bar\delta_i$ denote the numbers defined in Theorem \ref{many}.
It is clear that $d_i=\gcd(\bar\delta_0,\ldots,\bar\delta_i)
=\gcd(\delta_0,\ldots,\delta_i)$. We keep also this notation in the sequel.

\begin{lem}\label{cond}
The conductor of $S=\langle\bar\delta_0,\ldots,\bar\delta_k\rangle$ is
$$(\frac{d_0}{d_1}-1)(\bar\delta_1-d_1)+(\frac{d_1}{d_2}-1)(\bar\delta_2-d_2)+
\cdots+(\frac{d_{k-1}}{d_k}-1)(\bar\delta_k-d_k),$$ and $S$ is symmetric.
\end{lem}

\noindent{\bf Proof.} Since $\gcd(\bar\delta_0,\bar\delta_1)=d_1$,
we have that $i\bar\delta_1,0\le i\le\frac{d_0}{d_1}-1$, are all
different $\pmod{\bar\delta_0}$. They are also all smaller than
$\bar\delta_2$, since $\bar\delta_2>\frac{d_0}{d_1}\bar\delta_1$.
In the same way all $i\bar\delta_1+j\bar\delta_2,0\le i\le\frac{d_0}{d_1}-1,
0\le j\le\frac{d_1}{d_2}-1$ are all different $\pmod{\bar\delta_0}$, and
they are all smaller than $\bar\delta_3$, since $\bar\delta_3>
\frac{d_1}{d_2}\bar\delta_2>(\frac{d_1}{d_2}-1)\bar\delta_2+
(\frac{d_0}{d_1}-1)\bar\delta_1$ a.s.o. In this way we see that
the Apery set of $S$ with respect to $\bar\delta_0$
is $\{ j_1\bar\delta_1+j_2\bar\delta_2+\cdots+j_k\bar\delta_k;\, 0\le j_i<
\frac{d_{i-1}}{d_i},\, i=1,\ldots,k\}$ and $i_1\bar\delta_1+
i_2\bar\delta_2+\cdots+i_k\bar\delta_k>j_1\bar\delta_1+j_2\bar\delta_2+
\cdots+j_k\bar\delta_k$ if and only if $i_k=j_k,\ldots,i_s=j_s,i_{s-1}>
j_{s-1}$ for some $s$, i.e., if the last nonzero coordinate of
$(i_1-j_1,\ldots,i_k-j_k)$ is positive. (We have found
$\frac{d_0}{d_1}\frac{d_1}{d_2}\cdots\frac{d_{k-1}}{d_k}=\frac{d_0}{d_k}=
d_0=\bar\delta_0$ elements which are smallest in their congruence
classes $\pmod{\bar\delta_0}$.)
Hence, the largest number in the Apery set is
$\omega_{\bar\delta_0-1}=
(\frac{d_0}{d_1}-1)\bar\delta_1+(\frac{d_1}{d_2}-1)\bar\delta_2+
\cdots+(\frac{d_{k-1}}{d_k}-1)\bar\delta_k$. Since the conductor equals
$\omega_{\bar\delta_0-1}-(\bar\delta_0-1)$ (cf. Lemma \ref{Frobenius}),
 we get the first
statement after a small calculation. If $\omega_i=i_1\bar\delta_1+\cdots+
i_k\bar\delta_k$, it is easy to see that
$\omega_{\bar\delta_0-1-i}=(\frac{d_0}{d_1}-1-i_1)\bar\delta_1+\cdots+
(\frac{d_{k-1}}{d_k}-1-i_k)\bar\delta_k$. Thus $\omega_i+
\omega_{\bar\delta_0-1-i}=
\omega_{\bar\delta_0-1}$, which gives that $S$ is symmetric (cf. \cite{Ap}).

\medskip
For a semigroup $S$ and an integer $d>0$, we define the $d$-conductor
of $S$ to be $c_d(S)=\min\{ nd;\, md\in S\mbox{ if }m\ge n\}$. Thus $c_1(S)$
is the usual conductor of $S$.

\begin{cor}\label{dcond}
Let
$S=\langle\bar\delta_0,\ldots,\bar\delta_k\rangle$ and let $d_i=\gcd(
\bar\delta_0,\ldots,\bar\delta_i)$. Then
$$c_{d_i}(S)=(\frac{d_0}{d_1}-1)(\bar\delta_1-d_1)+
(\frac{d_1}{d_2}-1)(\bar\delta_2-d_2)+
\cdots+(\frac{d_{i-1}}{d_i}-1)(\bar\delta_i-d_i)$$
for every $i\le k$.
\end{cor}

\noindent{\bf Proof.} By the proof of Lemma \ref{cond}, the semigroup
$S_i=\langle\frac{\bar\delta_0}{d_i},\ldots,\frac{\bar\delta_i}{d_i}\rangle$
has conductor $c(S_i)=$
$$(\frac{d_0/d_i}{d_1/d_i}-1)(\frac{\bar\delta_1}{d_i}-\frac{d_1}{d_i})+
(\frac{d_1/d_i}{d_2/d_i}-1)(\frac{\bar\delta_2}{d_i}-\frac{d_2}{d_i})+
\cdots+(\frac{d_{i-1}/d_i}{d_i/d_i}-1)(\frac{\bar\delta_i}{d_i}-
\frac{d_i}{d_i}).$$
Then $c_{d_i}(\langle\bar\delta_0,\ldots,\bar\delta_k\rangle)=
d_ic(S_i)$. A calculation gives that
$\bar\delta_{i+1}>c_{d_i}(\langle\bar\delta_0,\ldots,\bar\delta_k\rangle)$,
hence $\bar\delta_j>c_{d_i}(\langle\bar\delta_0,\ldots,\bar\delta_k\rangle)$
if $j>i$. Thus $c_{d_i}(S)=
c_{d_i}(\langle\bar\delta_0,\ldots,\bar\delta_k\rangle)$.

\begin{lem}\label{conddelta}
For $i=2,\ldots,k$ we have $\bar\delta_i=\frac{1}{d_{i-1}}\sum_{j=1}^{i-1}
(d_{j-1}-d_j)\delta_j+\delta_i$. Thus the conductor of
$S=\langle\bar\delta_0,\ldots,\bar\delta_k\rangle$ is
$\sum_{i=1}^k(d_{i-1}-d_i)\delta_i+(1-d_0)$. Furthermore $c_{d_i}(S)=
\frac{1}{d_i}\sum_{j=1}^i\delta_j(d_{j-1}-d_j)+d_i-d_0$.
\end{lem}

\noindent{\bf Proof.} By a calculation, replacing in Lemma
 \ref{cond} and in Corollary \ref{dcond}  $\bar\delta_i$ with
$\frac{1}{d_{i-1}}\sum_{j=1}^{i-1}
(d_{j-1}-d_j)\delta_j+\delta_i$, we get the claim.

\medskip
For the next proposition, we
need a technical lemma. Let $g(t)=\sum_{i\ge0}a_it^i$, $a_0\ne0$
be a power series such that $\gcd(\{ i;\ a_i\ne0\})=1$. Let, for
$i=1,\ldots,k-1$, ${\bf d}_i=(d_i,\ldots,d_{k-1})$, and let
${\bf d}_i(g(t))=
(\epsilon_i(g),\ldots,\epsilon_{k-1}(g))$, where $\epsilon_s(g)=\min\{ j;\,
a_j\ne0,\mbox{ $d_s$ does not divide $j$}\}$. The easy proof of the
next lemma is left to the reader.

\begin{lem}\label{tech}
Let $g(t)=\sum_{i>0}a_it^i,a_0\ne0$, $h(t)=\sum_{i>0}b_it^i,b_0\ne0$, be
power series such that $\gcd(\{ i;\ a_i\ne0\}=\gcd(\{ i;\ b_i\ne0\}=1$.
Then\\
(a) ${\bf d}_i(gh)\ge\min({\bf d}_i(g),{\bf d}_i(h))$ (coefficientwise).\\
(b) If $g=h$ there is equality in (a).\\
(c) If ${\bf d}_i(g(t))=(\epsilon_1,\ldots,\epsilon_{k-1})$, then
${\bf d}_{i+1}((\sum_{i\ge\epsilon_1}a_it^i)/t^{\epsilon_1})=\\
(\epsilon_2(g)-\epsilon_1(g),\ldots,\epsilon_s(g)-\epsilon_1(g))$.
\end{lem}

We will call a power series {\em monic} if its least nonzero coefficient is 1.

\begin{prop}\label{ineq}
Let ${\mathcal O}=
{\mathbb C}[[t^{\delta_0}+\sum_{i\ge N}a_it^i,\sum_{i>\delta_0}b_it^i]]$
be a branch of standard parametrization and
with characteristic exponents $(\delta_0,\ldots,\delta_k)$. Let $\bar\delta_i$
be defined as in Theorem \ref{many}.
Then we have $\langle\bar\delta_0,\ldots,\bar\delta_k\rangle
\subseteq v({\mathcal O})$, i.e. $\bar\delta_i\in v({\mathcal O})$ for
$i=0,\ldots,k$.
\end{prop}

\noindent{\bf Proof.} Let, for $i=1,\ldots,k-1$,
${\bf d}_i=(d_i,\ldots,d_{k-1})$,
where $d_i=\gcd(\delta_0,\ldots,\delta_i)$ as above.
We will, by induction, construct monic elements $f_i
\in{\mathcal O}$ such that $v(f_i)=\bar\delta_i$ and such that
${\bf d}_i(f_i/t^{\bar\delta_i})=(\delta_{i+1}-\delta_i,
\delta_{i+2}-\delta_i,\ldots,\delta_k-\delta_i)$ if
$1\le i<k$. We let $f_0=t^{\delta_0}+
\sum_{i\ge N}a_it^i$. If $v(\sum_{i>\delta_0}b_it^i)$ is not a multiple
of $\delta_0$, then $o(\sum_{i>\delta_0}b_it^i)=\delta_1$
and we let $f_1=b_{\delta_1}^{-1}\sum_{i>\delta_0}b_it^i$.
If $v(\sum_{i>\delta_0}b_it^i)=m_0\delta_0$,
let $f_1'=\sum_{i>\delta_0}b_it^i-cf_0^{m_0}$,
where $c\ne0$ is chosen so that $v(f_1')>o(\sum_{i>\delta_0}b_it^i)$.
Repeat this until
$v(f_1^{(n)})=\delta_1$, and let $f_1=c'f_1^{(n)}$, where $c'$ is
chosen so that $f_1$ is monic. It is clear that
${\bf d}_1(f_1/t^{\bar\delta_1})=(\delta_2-\delta_1,
\delta_3-\delta_1,\ldots,\delta_k-\delta_1)$.
Suppose we have constructed
$f_0,f_1,\ldots,f_i\in{\mathcal O}$
so that the conditions in the proposition are
fulfilled. Then $f_i^{d_{i-1}/d_i}$ has value
$\gamma_i=\bar\delta_i\frac{d_{i-1}}{d_i}$, which is a multiple of $d_{i-1}$.
A simple calculation, using Lemma \ref{conddelta}, shows that
$\gamma_i-
c_{d_i}(\langle\bar\delta_0,\ldots,\bar\delta_i\rangle)=\delta_i-d_i+d_0>0$.
Thus of course $\gamma_i>
c_{d_{i-1}}(\langle\bar\delta_0,\ldots,\bar\delta_{i-1}
\rangle)$. This last means that $\gamma_i=\sum_{j=0}^{i-1} n_j\bar\delta_j$ for
some $n_j\ge0$. We choose $f_{i+1}'=f_i^{d_{i-1}/d_i}
-f_0^{n_0}\cdots f_{i-1}^{n_{i-1}}$.
From Lemma \ref{tech}(b) it follows
that ${\bf d}_i(f_i^{d_{i-1}/d_i}/t^{\bar\delta_i^{d_{i-1}/d_i}})=
{\bf d}_i(f_i/t^{\bar\delta_i})$. Since, for $j<i$,
${\bf d}_j(f_j/t^{\bar\delta_j})=(\delta_{j+1}-\delta_j,\ldots,\delta_k-
\delta_j)$, we have
${\bf d}_j(f_j/t^{\bar\delta_j})=(\delta_{i+1}-\delta_j,\ldots,\delta_k-
\delta_j)>(\delta_{i+1}-\delta_i,\ldots,\delta_k-\delta_i)$ (coefficientwise).
Lemma \ref{tech}(a) and (b) shows that
${\bf d}_i(f_0^{n_0}\cdots f_{i-1}^{n_{i-1}}/t^{\bar\delta_i})>
(\delta_{i+1}-\delta_i,\ldots,\delta_k-\delta_i)$.
Thus the smallest power in $f_{i+1}'$ which is not a multiple of $d_i$
and has nonzero coefficient is $\bar\delta_{i+1}$.
If $v(f_{i+1}')$ is not a multiple
of $d_i$, we choose $f_{i+1}=cf_{i+1}'$ ($c$ chosen so that $f_{i+1}$
is monic). If $v(f_{i+1}')$ is a multiple of $d_i$, then $\gamma_i>
c_{d_i}(\bar\delta_0,\ldots,\bar\delta_i\rangle)$ shows that $v(f_{i+1}'-
f_0^{m_0}\cdots f_i^{m_i})=v(f_{i+1}'')
>v(f_{i+1}')$ for some $m_0,\ldots,m_i\ge0$. We repeat until
$v(f_{i+1}^{(n)})=\bar\delta_{i+1}$, and let $f_{i+1}=c'f_{i+1}^{(n)}$,
where $c'$ is chosen so that $f_{i+1}$ is monic. It follows from
Lemma \ref{tech}(c) that
${\bf d}_{i+1}(f_{i+1}/t^{\bar\delta_{i+1}})=(\delta_{i+2}-\delta_{i+1},
\ldots,\delta_k-\delta_{i+1})$.

\begin{lem}\label{ocond} Let ${\mathcal O}$ be a branch with characteristic
exponents
$(\delta_0,\ldots,\delta_k)$.
Then the semigroup $v({\mathcal O})$ has conductor
$\sum_{i=1}^k(d_{i-1}-d_i)\delta_i+(1-d_0)$.
\end{lem}

\noindent{\bf Proof.} We make induction over the number $l$ of blowups we
need to get a regular branch. If $l=1$, then ${\mathcal O}=
{\mathbb C}
[[t^{\delta_0},t^{\delta_0+1}+\cdots]]$. It follows from Proposition \ref{ap}
that $v({\mathcal O})=\langle\delta_0,\delta_0+1\rangle$, which has
conductor $(\delta_0-1)\delta_0=
(\delta_0-1)(\delta_0+1)+1-\delta_0=(d_0-d_1)\delta_1+1-d_0$.
Suppose the claim is proved for $l-1$.
Let $c$ and $c'$ denote the
conductors of $v({\mathcal O})$ and $v({\mathcal O}')$, respectively.
In case a) of Lemma \ref{newpar},
a calculation using Lemma \ref{conddelta} gives $c-c'=\sum_{i=1}^k
(d_{i-1}-d_i)\delta_{i-1}=\delta_0^2-\delta_0$. By induction the
statement is true for $v({\mathcal O}')$. Proposition \ref{ap} shows it is
true for $v({\mathcal O})$. A similar calculation in case b) of Lemma
\ref{newpar} shows that
$c-c'=\delta_0^2-\delta_0$ also in this case. In case c)
of Lemma \ref{newpar} finally, we
get, by using $\delta_1-\delta_0=d_1, \delta_0=kd_1,\delta_1=(k+1)d_1$
for some $k$, that $c-c'=k^2d_1^2+kd_1=\delta_0^2-\delta_0$ also in
case c).

\medskip\noindent{\bf Proof of Theorem \ref{many}.} We know that
$\langle\bar\delta_0,\ldots,\bar\delta_k\rangle\subseteq v({\mathcal O})$
and that by Lemmas \ref{conddelta} and \ref{ocond}
these two semigroups have the same conductor. Since
$\langle\bar\delta_0,\ldots,\bar\delta_k\rangle$ is symmetric, all
strictly larger semigroups have smaller conductor.
This gives that the two semigroups are in fact the same.

\medskip
We get an easy criterion for a semigroup $\langle a_0,\ldots,a_k\rangle$
to be a semigroup for a plane branch. The following seems to be a simpler
characterization of the semigroup of a plane branch, with respect to
equivalent characterizations found in \cite{Za} or \cite{Br}.

\begin{prop}\label{semichar} Let $S$ be a semigroup which is minimally
generated by $a_0<a_1<\cdots<a_k$ and let $d_i=\gcd(a_0,\ldots,a_i),
i=0,\ldots,k$. Then
$S$ is the semigroup of a plane branch if and only if the following
conditions are satisfied.\\
(a) $d_0>d_1>\cdots>d_k=1$.\\
(b) $a_i>{\rm lcm}(d_{i-2},a_{i-1})$ for $i=2,\ldots,k$.
\end{prop}

\noindent{\bf Proof.} The necessity follows from Theorem \ref{many}, the
sufficiency from the branch ${\mathbb C}[[t^{a_0},t^{a_1}+
t^{a_1+a_2-{\rm lcm}(d_0,a_1)}+
\cdots+t^{a_1+\cdots+a_k-({\rm lcm}(d_0,a_1)+\cdots+
{\rm lcm}(d_{k-2},a_{k-1}))}]]$.

\medskip
We give two concrete examples.

\medskip
\noindent{\bf Example} Let $S=\langle 30,42,280,855\rangle$.
Then $S$ satisfies the conditions in Proposition \ref{semichar},
so $S=v({\mathcal O})$ for some ${\mathcal O}$. We can choose e.g.
${\mathcal O}={\mathbb c}[[t^{30},t^{42}+t^{112}+t^{127}]]$. The
conductor equals $t^{1554}{\mathbb C}[[t]]$. With the notation of the
previous section, the
multiplicity sequence is  $M(30,42),M(6,70),M(2,15),\dots$, which is
$30,12^{(2)},6^{(13)},4,2^{(9)},1^{(2)},\ldots$.

\medskip
\noindent{\bf Example} Let
$\mathcal O={\mathbb C}[[x,y]]$ with
$$x=t^{2\cdot 2^n},\, y=t^{3\cdot 2^n}+t^{3\cdot 2^n+2^{n-1}}+\cdots+
t^{3\cdot 2^n+2^{n-1}+\cdots+2+1}.$$
The
generators of $v({\mathcal O})$ are $\bar\delta_0=2^{n+1},\bar\delta_i=
2^{n-i+1}(3\cdot2^{2i-2}+(4^{i-1}-1)/3)$ for $i=1,\ldots,n+1$.

\section{COMPLETE INTERSECTION RINGS\\ ARISING FROM THE SEMIGROUP OF A PLANE 
BRANCH}

Let $S=\langle \bar\delta_0,\ldots,\bar\delta_k\rangle=v({\mathcal O})$
be the semigroup of a plane branch, where
$\bar\delta_0<\bar\delta_1<\ldots<\bar\delta_k$ is a minimal set of
generators of
$S$, and let
${\mathbb C}[S]={\mathbb C}[t^{\bar\delta_0},\ldots,t^{\bar\delta_k}]=
{\mathbb C}[Y_0,\ldots,Y_k]/I=T$. We will show that $T$
has an associated graded ring (in the $(Y_0,\ldots,Y_k)$-filtration),
which is a complete intersection.  In particular this implies that
$T$ is a complete intersection \cite{VV}. We will use
\cite[Theorem 1]{Ro}
which states that if all elements in ${\rm Ap}(S,\bar\delta_0)$,
the Apery set of $S$ with respect to $\bar\delta_0$, have unique
expressions as linear combinations of the generators of $S$, then the
relations are determined by the minimal elements above the Apery set.
In the following results, we suppose $S=v({\mathcal O})$, where ${\mathcal
O}$ is
a plane branch. We also keep the notation of the previous sections.

\begin{lem}  All elements
in ${\rm Ap}(S,\bar\delta_0)$ have unique expressions.
\end{lem}

\noindent{\bf Proof.} The elements in ${\rm Ap}(S,\bar\delta_0)$ are of the
form $i_1\bar\delta_1+\cdots+i_k\bar\delta_k$, with $0\le i_j<d_{j-1}/d_j$
(cf. proof of Lemma \ref{cond}).
Suppose $i_1\bar\delta_1+\cdots+i_k\bar\delta_k=
j_0\bar\delta_0+\cdots+j_k\bar\delta_k$. Then $i_k\bar\delta_k\equiv
j_k\bar\delta_k\pmod{d_{k-1}}$. Since
$i_1\bar\delta_1+\cdots+i_{k-1}\bar\delta_{k-1}<\bar\delta_k$, this
implies that $i_k=j_k$. If $k>1$ we get $i_{k-1}\bar\delta_{k-1}=
j_{k-1}\bar\delta_{k-1}\pmod{d_{k-2}}$, which gives $i_{k-1}=j_{k-1}$
a.s.o. Finally $0=j_0\bar\delta_0$, so $j_0=0$.

\medskip
Next we determine the ``minimals" (cf. \cite{Ro}), i.e. the minimal
elements $(n_1, \cdots,n_k) \in {\mathbb N}^k$ such that
$ n_1\bar\delta_1+\cdots+n_k\bar\delta_k\notin
{\rm Ap}(S,\bar\delta_0)$ (the order in ${\mathbb N}^k$ is the usual one).
Some
$n_j$ must be at least $d_{j-1}/d_j$,
otherwise the element belongs to ${\rm Ap}(S,\bar\delta_0)$. On the other
hand at most one $n_j\ge d_{j-1}/d_j$ and there must be equality, if
the element is minimal outside ${\rm Ap}(S,\bar\delta_0)$. Thus the
minimals are
$$\{(d_0/d_1,0,\cdots,O),(0,d_1/d_2,0,\cdots,0), \cdots,
(0,\cdots,0,d_{k-1}/d_k)\}.$$
Thus the following theorem follows from
\cite[Theorem 1]{Ro}.

\begin{thm} A minimal presentation for ${\mathbb C}[S]$ is
  $${\mathbb C}[S]={\mathbb C}[Y_0,\ldots,Y_k]/(Y_1^{d_0/d_1}
-m_1,\ldots,Y_k^{d_{k-1}/d_k}-m_k)$$
where $m_j$ is a monomial in $Y_0,\ldots,Y_j$ for $j=1,\ldots,k$. Thus
${\mathbb C}[S]$ is a complete intersection.
\end{thm}

\begin{cor} The associated graded ring of ${\mathbb C}[S]$ with respect to
the filtration given by powers of $(Y_0,\ldots,Y_k)$ is
${\mathbb C}[Y_0,\ldots,Y_k]/(Y_1^{d_0/d_1},\ldots,Y_k^{d_{k-1}/d_k})$. Thus
it is a complete intersection.
\end{cor}

\noindent{\bf Proof.} Since $m_j=Y_0^{n_0}\cdots Y_{j-1}^{n_{j-1}}$ and
$n_0\bar\delta_0+\cdots+n_{j-1}\bar\delta_{j-1}=(d_{j-1}/d_j)\bar\delta_j$, 
it is clear that $n_0+\cdots+n_{j-1}>(d_{j-1}/d_j)$, so
${\rm in}(Y_j^{d_{j-1}/d_j}-m_j)=Y_j^{d_{j-1}/d_j}$. Since $Y_1^{d_0/d_1},
\ldots,Y_k^{d_{k-1}/d_k}$ is a regular sequence, we get the result,
cf. \cite{VV}.

\medskip
\noindent{\bf Remark.} Notice that not only for semigroups of plane
branches the two results above hold. For example, if $S=\langle
4,6,7\rangle$, then
$S$ is not the semigroup of a plane branch, but ${\mathbb C}[S] = {\mathbb
C} [X,Y,Z]/(Y^2 - X^3, Z^2 - X^2Y)$ is a complete intersection and also
its associated graded ring is a complete intersection.

\begin{cor} The generating function for $S$, i.e. $\sum_{i\in S}t^i$,
equals
$$(1-t^{(d_0/d_1)\bar
\delta_1})\cdots(1-t^{(d_{k-1}/d_k) \bar
\delta_k})/((1-t^{\bar\delta_0})\cdots (1-t^{\bar\delta_k})).$$
\end{cor}

\noindent{\bf Proof.} As graded algebra ${\mathbb C}[S]$ is generated
by $k+1$ elements of degrees $\bar\delta_i$, $i=0,\ldots,k$ and has $k$
minimal relations of degrees $(d_{i-1}/d_i)\bar\delta_i$, $i=1,\ldots,k$,
which constitute a regular sequence.

\medskip
\noindent{\bf Examples.} If ${\mathcal O}={\mathbb C}[[t^8,t^{12}+t^{14}+
t^{15}]]$, then $v({\mathcal O})=\langle8,12,26,53\rangle$ so the generating
function is $(1-t^{24})(1-t^{52})(1-t^{106})/((1-t^8)(1-t^{12})(1-t^{26})
(1-t^{53}))$.

\smallskip
If ${\mathcal O}={\mathbb C}[[t^{30},t^{42}+t^{112}+
t^{127}]]$, then $v({\mathcal O})=\langle30,42,280,855\rangle$
so the generating
function is
$(1-t^{210})(1-t^{840})(1-t^{1710})/((1-t^{30})(1-t^{42})(1-t^{280})
(1-t^{855}))$.


\begin{thebibliography}{Dillo 83}

\bibitem[1]{Sa} Samuel, P. {Alg\'ebricit\'e de certains points
singuliers alg\'ebro\"{\i}d}. J. Math. Pures Appl. {\bf 1956}, {\em 35}, 1--6.

\bibitem[2]{Za} Zariski, O. {Le probl\`eme des modules pour les
branches planes}. Hermann, Paris, 1986.

\bibitem[3]{Ap} Ap\'ery, R. {Sur les branches superlin\'eaires
des courbes alg\'ebriques}. C. R. Acad. Sci. Paris {\bf 1946}, {\em 222},
1198--2000.

\bibitem[4]{Hi} Hironaka, H. {On the arithmetic and effective
genera of algebraic curves}.  Mem. Coll. Sci. Univ. Kyoto {\bf 1957}, 
{\em 30}, 177--195.

\bibitem[5]{dV} du Val, P. {The Jacobian algorithm and the
multiplicity sequence of an algebraic branch}. Revue Fac. Sci. Univ.
Istanbul, Ser. A, {\bf 1942}, {\em tome VII, fasc. 3--4}, 107--112.

\bibitem[6]{Ca} Campillo, A. {Algebroid curves in positive
characteristic}. Lect. Notes in Math. {\bf 1980}, {\em 813}, Springer.

\bibitem[7]{Ab} Abhyankar, S. {Desingularization of plane curves}.
Proc. Sympos. Pure Math. {\bf 1983}, {\em 40}, 1--45.

\bibitem[8]{Be-Ca} Bertin, J.; Carbonne, P. {Semi-groupes
d'entiers et application aux branches}. J. Algebra {\bf 1977}, 
{\em 49}, 81--95.

\bibitem[9]{Me} Merle, M. {Invariants polaires des courbes planes}.
Invent. Math. {\bf 1977}, {\em 41}, 103--111.

\bibitem[10]{Br} Bresinsky, H.  {Semigroups corresponding to
algebroid branches in the plane}. Proc. Am. Math. Soc. {\bf 1972}, {\em 32},
381--384.

\bibitem[11]{VV} Valabrega, P.; Valla, G. {Form rings and regular
sequences}. Nagoya Math. J. {\bf 1978}, {\em 72}, 93--101.

\bibitem[12]{Ro} Rosales, J.-C. {Numerical semigroups with Ap\'ery
sets of unique expression}. J. Algebra {\bf 2000}, {\em 226}, 479--487.

\end{thebibliography}
\end{document}